\documentclass{amsart}
\usepackage{ae,aecompl}\usepackage{ifthen}
\usepackage[utf8]{inputenc}
\usepackage[a4paper,twoside]{geometry}
\usepackage[alphabetic]{amsrefs}\usepackage{pifont}
\usepackage{amsthm}
\usepackage{graphicx}
\usepackage{amssymb}
\usepackage{url}
\newboolean{dum}
\setboolean{dum}{false}
\newcommand{\comment}[1]{\ifthenelse{\boolean{dum}}{
{\par\noindent\Huge\ding{46}} \fbox{\parbox{10cm}{#1}}\par}{}}
\newcommand{\clabel}[1]{\comment{label #1}}
\newtheorem{lemma}[subsection]{Lemma}
\newtheorem{thm}[subsection]{Theorem}
\newtheorem{prop}[subsection]{Proposition}
\newtheorem{cor}[subsection]{Corollary}
\numberwithin{equation}{section}
\newcommand{\cM}{\mathcal{M}}
\newcommand{\cX}{\mathcal{X}}
\DeclareMathOperator{\vol}{vol}
\DeclareMathOperator{\card}{card}
\DeclareMathOperator{\sgn}{sgn}
\DeclareMathOperator{\Tr}{Tr}
\DeclareMathOperator{\Ai}{Ai}
\DeclareMathOperator{\Real}{Re}
\DeclareMathOperator{\linspan}{span}
\DeclareMathOperator{\Exp}{exp}
\begin{document}
\comment{
$ $Id: article.tex,v 1.21 2008/04/21 12:02:03 enord Exp $ $
}
\title{The Anti-Symmetric GUE Minor Process}
\author{Peter J. Forrester \and Eric Nordenstam}
\address{Department of Mathematics and Statistics \\
The University of Melbourne 
Victoria 3010, Australia}
\maketitle

\begin{quote}
\small{\sc Abstract.}
Our study is initiated by a multi-component particle system underlying
the tiling of a half hexagon by three species of rhombi. In this particle
system species $j$ consists of $\lfloor j/2 \rfloor$ particles which are
interlaced with neigbouring species. The joint probability density
function (PDF) for this particle system is obtained, and is shown in a
suitable scaling limit to coincide with the joint eigenvalue PDF
for the process formed by the successive minors of anti-symmetric GUE
matrices, which in turn we compute from first principles. The correlations
for this process are determinantal and we give an explicit formula for the
corresponding correlation kernel in terms of Hermite polynomials.
Scaling limits of the latter are computed, giving rise to the
Airy kernel, extended Airy kernel and bead kernel at the soft edge and in
the bulk, as well as a new kernel at the hard edge.
\end{quote}

\section{Introduction}
A classical problem in random matrix theory is to compute the statistical properties of
the eigenvalues, given the distribution of the entries of the matrix. For the well
known Gaussian orthogonal, unitary and symplectic ensembles (GOE, GUE and GSE; see
e.g.~\cite{forrester:book}) this problem can be solved exactly, with $k \times k$
determinant, and $2k \times 2k$ Pfaffian formulas available for the $k$-point
correlations in the case of the GUE, and GOE, GSE respectively. Of these the
GUE is special because determinantal expressions can also be derived for the
correlations of the process specified by the eigenvalues of all the successive
minors \cites{johansson:gue, okounkov:birth}, whereas for the GOE and GSE the
correlations for the minors are not known.

While interest in the eigenvalues of GUE matrices stems from its use in modelling
the highly excited energy levels of classically chaotic quantum systems with no time
reversal symmetry, the interest in the GUE minor process stems from a work of
Baryshnikov \cite{baryshnikov:gue}
relating to queues. Subsequently the GUE minor process has 
been identified by Johansson and Nordenstam within
a number of random tiling problems \cite{johansson:gue}.
In the language of stepped surfaces (which are equivalent
to certain tilings) similar applications were noticed in~\cite{okounkov:birth},
while in \cite{borodin:png_tasep}
the GUE minor process is encountered in the asymmetric
exclusion process.

One of the random tiling problems identified in \cite{johansson:gue}
as relating to the GUE minor process is the tiling of a hexagon by three species of
rhombi. This tiling problem is equivalent to non-intersecting lattice paths which
at each step move up one unit or down one unit in a so called watermelon formation
(see \cite{krattenthaler:walkers}). A very natural variation of the non-intersecting
paths model is to impose the boundary condition that the paths cannot go below
the $x$-axis \cites{forrester89,krattenthaler:walkers,gillet}, 
which in fact corresponds to the
tiling of a half hexagon by the same three species of rhombi as referred
to above (see Figure \ref{fig:hexagon} below).

This variation initiates the study of the present work, where we begin in
section 2 by analyzing such non-intersecting paths from the
same viewpoint which in the study \cite{johansson:gue} gave rise to the GUE minor
process. We obtain in corollary \ref{thm:AS} a joint probability 
density function (PDF)  which can be interpreted as a multi-component 
particle system in which the number of particles in species $j$ consists of
$\lfloor j/2 \rfloor$ particles $(j=1,\dots,2n+1)$. The marginal distribution
of species $j$ can also be computed exactly
(see theorem \ref{As}), and is precisely that of a $j \times j$ anti-symmetric GUE matrix
(see e.g.~\cite{forrester:book}*{Ch.~1}). By analogy with the findings from
\cite{johansson:gue}, this motivates us to seek the joint eigenvalue PDF of the
anti-symmetric GUE minor process, and to compare this against the
joint PDF obtained from analyzing the non-intersecting lattice paths
model. We do this in Section 3, and find the two joint PDFs are
identical.

It should be noted that the discrete processes that in this article and 
in \cite{johansson:gue} give rise to random matrix distributions is not 
the positions of the random walkers, but the holes between them, se section 2. 
It is also possible to obtain random matrix distributions by looking at
the positions of random walkers directly, see for examble \cites{gillet,kuijlaars}
for examples of such models. 
For tilings of a hexagon, it is well known that the limit
shape of the disordered region is a circle,~\cite{cohn:shape_typical}.
We believe that the limit shape in our half-hexagon  model is a semicircle, 
but do not prove it.

As already mentioned, the GUE minor process is determinantal, with the correlations
between eigenvalues from prescribed minors being given by an $r \times r$ determinant.
Explicitly, with the coordinate $(s,y)$ denoting an eigenvalue from the $s$-th principal
minor having value $y$, it was found that \cite{johansson:gue} 
\begin{equation}
\rho_{(r)} (\{ (s_j,y_j) \}_{j=1, \dots, r} )=
\det[K^{\rm GUEm}((s_j, y_j), (s_k, y_k)) ]_{j,k=1,\dots, r}
\end{equation}
with
\clabel{eq:58}\begin{equation}
\label{eq:58}
K^{\rm GUEm}((s, x), (t, y)) =
\begin{cases}
\displaystyle
\frac{e^{-(x^2+y^2)/2} }{\sqrt{\pi}}
\sum_{k=1}^t
\frac{H_{s-k}(x)H_{t-k}(y)}{2^{t-k} (t-k)!},  &s\geq t\\
\quad\\
\displaystyle
-\frac{e^{-(x^2+y^2)/2} }{\sqrt{\pi}}
\sum_{k=-\infty}^0
\frac{H_{s-k}(x)H_{t-k}(y)}{2^{t-k} (t-k)!},  &s\geq t
\end{cases}
\end{equation}
(here $H_j(x)$ denotes the Hermite polynomials of degree $j$, satisfying
$\int_{\mathbb{R}} H_i(x) H_j(x) e^{-x^2}\, dx = \delta_{ij} n!2^n\sqrt\pi$).
In section 4, proposition \ref{thm:cont_borodin}, we show that similarly
for the anti-symmetric GUE minor process
\begin{equation}\label{eq58a}
\rho_{(r)} (\{ (s_j,y_j) \}_{j=1, \dots, r} )=
\det[K^{\rm aGUEm}((s_j, y_j), (s_k, y_k)) ]_{j,k=1,\dots, r}
\end{equation}
with
\begin{equation}
K^{\rm aGUEm}((s, x), (t, y)) = 
\begin{cases}
\displaystyle
\frac{2 e^{-(x^2+y^2)/2}}{\sqrt{\pi}} 
\sum_{l=1}^{\lfloor t/2\rfloor} \frac{H_{s - 2l}(x) H_{t - 2l}(y)}{
2^{t - 2l} (t - 2l)! },& s \ge t, \\
\quad\\
\displaystyle 
- \frac{2 e^{-(x^2+y^2)/2}}{\sqrt{\pi}} 
\sum_{l=-\infty}^{0} \frac{H_{s - 2l}(x) H_{t - 2l}(y)}{
2^{t - 2l} (t - 2l)! },& s < t.
\end{cases}
\end{equation}

In section 5 we turn our attention to various scaled limits of
(\ref{eq58a}). In particular, two distinct soft edge scaling
limits (neighbourhood of the largest eigenvalues of the respective
minors), as well as a bulk and hard edge scaling (the latter referring
to the eigenvalues in the neighbourhood of the origin) are computed.
At the soft edge, 
with the species (i.e.~minors) labels chosen to differ by $2n$ by a fixed amount,
the limiting correlation kernel is found to be the well known
Airy kernel (see e.g.~\cite{forrester:book}*{Ch.~4})
\begin{equation}\label{Ksoft}
K^\text{soft}(x,y) =
\int_0^\infty \Ai(x+u) \Ai(y+u)\, du
\end{equation}
independent of the species. The other soft edge scaling limit
analyzed is when the species differ by $O(n^{2/3})$. In the case
of the GUE minor process,
study of this limit is motivated by a mapping
to the interface of a certain droplet polynuclear growth model
\cite{prahofer:png_droplet}. Fluctuations of the latter are governed
by the dynamical extension of the Airy kernel
\begin{equation}
\label{eq:34}
K^\text{soft}((\tau_x, x), (\tau_y, y)) = 
\begin{cases}
\int_0^\infty e^{-(\tau_y-\tau_x)u }
\Ai(x+u) \Ai(y+u) \, du, &\tau_y\geq \tau_x\\
-\int_{-\infty}^0 e^{-(\tau_y-\tau_x)u }
\Ai(x+u) \Ai(y+u) \, du, &\tau_y< \tau_x.
\end{cases}
\end{equation}
It is (\ref{eq:34}) which again appears in this second soft edge scaling.
In the bulk scaling of the GUE minor process, the limiting
correlation kernel was found to be \cite{forrester:nagao}
\begin{equation}
\label{eq:33}
K^{\rm bead}( (\tau_x, x), (\tau_y, y)) = 
\begin{cases}
\displaystyle
\frac{1}{2}\int_{-1}^1 (is)^{\tau_y-\tau_x}
 e^{is\pi (s-y)}\,ds,&\tau_y>\tau_x\\
\displaystyle
-\frac{1}{2}\int_{\mathbb{R}\setminus [-1,1]} (is)^{\tau_y-\tau_x} 
e^{is\pi (s-y)}\,ds,&\tau_y<\tau_x.\\
\end{cases}
\end{equation}
This corresponds to the $\gamma=0$ (isotropic) case of the
correlation kernel for the bead process~\cite{boutillier}.
It too is reclaimed in the appropriate bulk scaling limit
of the anti-symmetric GUE minor process. The hard edge scaling
limit has no analogue in the GUE minor process. The limiting
correlations are given in proposition \ref{prop5.5}.

\smallskip
\noindent
{\it Note:} 
While this work was being prepared, a research
announcement by M. Defosseaux was posted on the arXiv, \cite{defosseaux},
stating  results equivalent to our theorem~\ref{thm:aGUE} 
and proposition~\ref{thm:cont_borodin}.

\section{The half hexagon}
\begin{figure}
\includegraphics{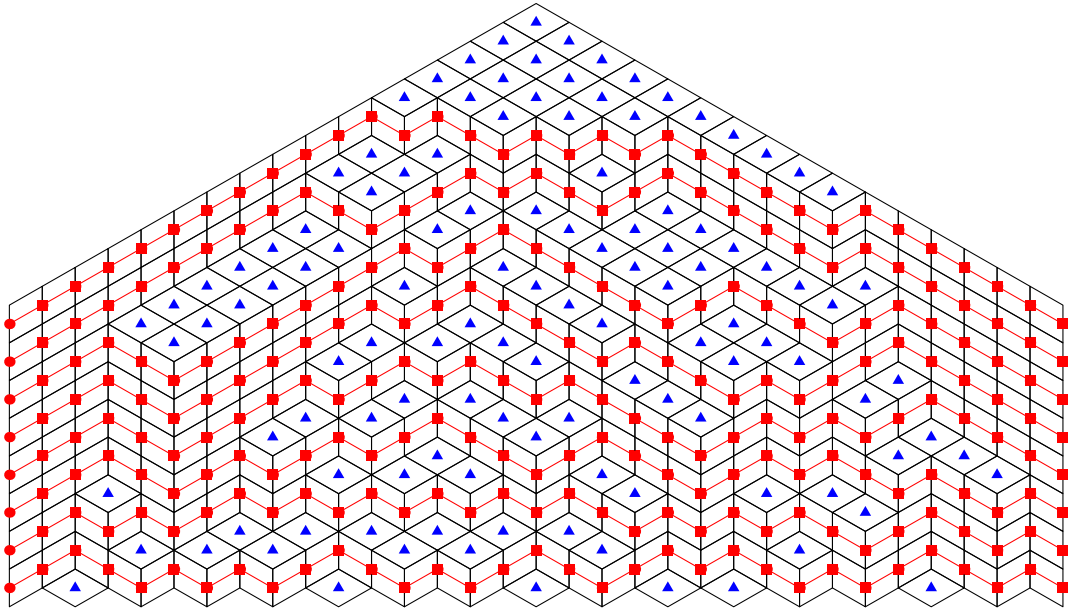}
\caption{Tiling with rhombuses of ``half a hexagon''
and the corresponding non-intersecting random walks. 
The positions of blue particles are marked by triangles 
and the red particles by squares. 
The picture was generated with the method of~\cite{propp}. 
 }
\label{fig:hexagon}
\end{figure}
Consider $p$ simple symmetric random walks started at
$(0,2i-2)$
conditioned to end up at $(2N,2i-2)$, for $i=1$, \dots, $p$.
They are conditioned never to intersect and  never to go below 
the $x$ axis.
Such configurations are in bijection to tilings of a half-hexagon, as
seen in the Figure \ref{fig:hexagon}. 
To count such configurations we need the following result
of Krattenthaler et al.~\cite{krattenthaler:walkers}.

\begin{prop}[Theorem 6 in \cite{krattenthaler:walkers}]
Let $e_1<e_2<\cdots<e_p$ with $e_i= m \pmod{2}$, $i=1,2,\dots, p$.
The number of stars with $p$ branches, the $i$-th branch running
from $A_i=(0,2i-2)$ to $E_i=(m,e_i)$, $i=1,2,\dots, p$,
and never going below the $x$-axis equals
\begin{equation}
\prod_{i=1}^p
\frac{(e_i+1)(m+2i-2)!}{
(\frac{m+e_i}{2} +p)!(\frac{m-e_i}{2} +p-1)!} 
\prod_{1\leq i<j\leq p}
\left(
\left(\frac{e_i+1}{2}\right)^2-
\left(\frac{e_j+1}{2}\right)^2
\right).
\end{equation}
\end{prop}
Here, a star means a configuration of $p$ simple symmetric random
walks, or branches.
Given that the red particles on line $x=n$ occur at
positions $e_1<e_2<\cdots <e_p$, the number of ways 
to tile the area to the left of line $n$ is
\begin{equation}
\prod_{i=1}^p
\frac{(e_i+1)(n+2i-2)!}{
(\frac{n+e_i}{2} +p)!(\frac{n-e_i}{2} +p-1)!} 
\prod_{1\leq i<j\leq p}
\left(
\left(\frac{e_i+1}{2}\right)^2-
\left(\frac{e_j+1}{2}\right)^2
\right).
\end{equation}
The number of ways to tile the area to the right of line $n$
is 
\begin{equation}
\prod_{i=1}^p
\frac{(e_i+1)(2N-n+2i-2)!}{
(\frac{2N-n+e_i}{2} +p)!(\frac{2N-n-e_i}{2} +p-1)!} 
\prod_{1\leq i<j\leq p}
\left(
\left(\frac{e_i+1}{2}\right)^2-
\left(\frac{e_j+1}{2}\right)^2
\right).
\end{equation}

Thus the probability of a certain configuration of red
particles on line $n$ is the product of the above two expressions 
divided by the total number of tilings of the area. 

\comment{label thm:red}
\begin{lemma}
\label{thm:red}
The probability distribution of red dots along line $x=n$
is
\begin{multline}
\mathcal{P}^n_{\rm red} [e_1<e_2<\dots<e_p] =\\ 
Z^{-1} \prod_{1\leq i<j\leq p} ( (e_j+1)^2-(e_i+1)^2 )^2
\prod_{i=1}^p
\frac{(e_i+1)^2(n+2i-2)!(2N-n+2i-2)!}{
(\frac{n+e_i}{2}+p)!(\frac{2N-n+e_i}{2} + p)! (\frac{n-e_i}{2}+p-1)!
(\frac{2N-n-e_i}{2} + p-1)!}
\end{multline} 
for $0\leq e_1< \dots < e_p \leq 2p+n-2$ and $e_i=n\bmod{2}$ for $i=1,\dots,p$.
\end{lemma}

We are more interested in the distribution of the holes, i.e. the blue dots. 
This can be computed using a result from \cite{borodin:duality}.
To set this up, 
let $X=\{x_1,\dots, x_M\}$ be a set of $M$ real numbers. 
Consider measures on the following form.
For $w:X\rightarrow\mathbb{R}$,
\begin{equation}
\mathcal{P}_w^{(m)}[x_{i_1}, \dots, x_{i_m}] = 
\text{const} \prod_{1\leq k<l\leq m} (x_{i_k} - x_{i_l})^2 
\prod_{k=1}^m w(x_{i_k})
\end{equation}
is a probability measure which is supported on sets of $m$ points. 
Notice that the measure above for the red particles is on this form.
Denote the measure on the holes by 
\begin{equation}
\bar{\mathcal{P}}_w^{M-m}(A) = \mathcal{P}_w^{m}(X\setminus A),
\end{equation}
this is of course a measure supported on sets of $M-m$ elements.
Borodin proves the following formula.

\comment{label eqn:borodin\_duality}
\begin{prop}[Proposition 2 in \cite{borodin:duality}]
\label{eqn:borodin_duality}
Let $u, v:X\rightarrow\mathbb{R}$ be functions such that 
for all $x_k\in X$, 
\begin{equation}
u(x_k) v(x_k) =\frac{1}{ \prod_{i\neq k} (x_k-x_i)^2}.
\end{equation}
Then
\begin{equation}
\bar{\mathcal{P}}_u^{M-m}=
\mathcal{P}_v^{M-m}.
\end{equation}
\end{prop}

With this tool we can now compute the distribution of blue particles.

\comment{ label thm:blue}
\begin{thm}
\label{thm:blue}
On line $x=n\leq N$, for $n$ even, the distribution of the 
blue particles is 
\begin{multline}
\mathcal{P}^n_{\rm blue}\{ e_1,e_2,\dots,e_{p'}\} =\\ 
\tilde Z^{-1} \prod_{1\leq i<j\leq p'} ( (e_j+1)^2-(e_i+1)^2 )^2
\prod_{i=1}^{p'}
\frac{(\frac{2N-n+e_i}{2} + p)! 
(\frac{2N-n-e_i}{2} + p-1)!}{
(n+2i-2)!(2N-n+2i-2)!(\frac{n+e_i}{2}+p)!(\frac{n-e_i}{2}+p-1)!}
\end{multline}
for $p'=n/2$ and $0\leq e_1<e_2<\dots<e_{p'}\leq 2p+n-2$ 
and $e_i$ even for all $i$.
For odd $n$ it is 
\begin{multline}
\mathcal{P}^n_{\rm blue}\{ e_1,e_2,\dots,e_{p'}\} =\\ 
\tilde Z^{-1} \prod_{1\leq i<j\leq p'} ( (e_j+1)^2-(e_i+1)^2 )^2
\prod_{i=1}^{p'}
\frac{(\frac{2N-n+e_i}{2} + p)! 
(\frac{2N-n-e_i}{2} + p-1)!(e_i+1)^2}{
(n+2i-2)!(2N-n+2i-2)!(\frac{n+e_i}{2}+p)!(\frac{n-e_i}{2}+p-1)!}
\end{multline}
for $p' = (n-1)/2$ and again $1\leq e_1<e_2<\dots<e_{p'}\leq 2p+n-2$
and all $e_i$ odd for all $i$.
\end{thm}
\begin{proof}
For $n$ even, let $X=\{0, 2, \dots, 2M-2\}$. 
By proposition \ref{eqn:borodin_duality}, we need to compute
\begin{equation}
\prod_{\begin{smallmatrix}
y\in X\\ y\neq x
\end{smallmatrix}} 
 ((x+1)^2-(y+1)^2 )^{-2} =
2^{-4M}
\left( 
\frac{(M-\frac{x}{2})! (M+1 +\frac{x}{2})! }{
x+1}
\right)^{-2}.
\end{equation}
For $n$ odd,  $X=\{1, 3, \dots, 2M-1\}$ and we need to compute
\begin{equation}
\prod_{\begin{smallmatrix}
y\in X\\ y\neq x
\end{smallmatrix}} 
 ((x+1)^2-(y+1)^2 )^{-2} =
2^{-4M}
\left( 
\frac{(M-\frac{1+x}{2})! (M+\frac{1+x}{2})! }{
(x+1)^2}
\right)^{-2}. 
\end{equation}
Combining these with the formula in 
lemma~\ref{thm:red} gives the result.
\end{proof}

Let us fix some notation.
For any $k\leq n$, let the positions of the $\lfloor k/2\rfloor$ 
blue particles on line $k$ 
be denoted $x^{(k)}_{1}>\dots>x^{(k)}_{\lfloor k/2\rfloor}$.
Let $x^{(k)}=(x^{(k)}_{1},\dots,x^{(k)}_{\lfloor k/2\rfloor})$.
For simplicity of notation later, let $x^{(2l+1)}_{l}=0$ 
for  $l=0$, 1,\dots, $\lfloor n/2\rfloor$. 
Obviously, for $k$ odd, $x^{(k)}_i$ is always even and for $k$ even, 
$x^{(k)}_i$ is always odd.

It so happens, as is easily seen from the picture, 
that the blue particles  must 
fulfill certain interlacing requirements,
\begin{equation}
x^{(k)}_{i+1}<x^{(k-1)}_i<x^{(k)}_i
\end{equation}
for all $k=2$, \dots, $n$ and $i=1$, \dots, $\lfloor k/2\rfloor$. 
Henceforth we will write $x^{(k+1)}\succ x^{(k)}$ when $x^{(k+1)}$ and $x^{(k)}$
satisfy the above interlacing. 
Let $K(x^n)=\{(x^{(1)}, \dots, x^{(n)}):x^{(n)}\succ x^{(n-1)}\succ\dots
\succ x^{(1)}\}$.

The next observation to make is that the distribution $x^{(1)}$, 
\dots, $x^{(n)}$, 
given $x^n$ is uniform in the cone given by the above inequalities. 
\begin{prop}
The joint probability of all blue particles on the first $n$ lines is
\begin{equation}
\mathcal{P}_{\rm blue}^{(1,n)} \{ x^{(1)},\dots, x^{(n)}\} = 
\frac{\chi_1(x^{(1)}, x^{(2)})
\dots
\chi_{n-1}(x^{(n-1)}, x^{(n)})}{
\card K(x^{(n)})}
\mathcal{P}^n_{\rm blue} \{
x^{(n)}
\}
\end{equation}
where 
$\chi_k(x^{(k)}, x^{(k+1)})=\mathbf{1}\{ x^{(k+1)}\succ x^{(k) }\}$
and $\card K(x^{(n)})$ is the cardinality of $K(x^{(n)})$.
\end{prop}

Seeing this, one would immediately want to compute $\card K(x^{(n)})$. 

\begin{lemma}\label{lem2.6}
Let $Z_2=1$, $Z_n = \prod_{j=1}^{n-2} j!!$, $n \ge 3$.
For $ n$ even, 
\begin{equation}\label{2.15}
\card K(x^{(n)}) = Z_n^{-1}\prod_{1\leq i<j\leq n/2} 
\left(
\frac{(x^{(n)}_i)^2}{4}- \frac{(x^{(n)}_j)^2}{4}
\right)
\end{equation}
and for $ n$ odd ($n \ge 3$), 
\begin{equation}\label{2.16}
\card K(x^{(n)}) = Z_n^{-1} \prod_{1\leq i<j\leq (n-1)/2} 
\left(
\frac{(x^{(n)}_i)^2}{4}- \frac{(x^{(n)}_j)^2}{4}
\right)
\prod_{i=1}^{(n-1)/2} \frac{x^{(n)}_i}{2}
\end{equation}
\end{lemma} 
\begin{proof}
The proof is an inductive one. With the products over $i < j$ interpreted as unity
in the cases $n=2$ and 3, we see by inspection that both (\ref{2.15}) and
(\ref{2.16}) are correct in these base cases.

Suppose (\ref{2.15}) has been established for $n=2N$. Then we see that
\begin{equation}\label{2.15a}
\card K(x^{(2N+1)}) = \sum_{x^{(2N)}: x^{(2N+1)} \succ x^{(2N)}} \card
K(x^{(2N)}).
\end{equation}
Use of the Vandermonde determinant identity in (\ref{2.15}) shows
$$
\card
K(x^{(2N)}) = Z_{2N}^{-1} \det \left [ \left (
\frac{(x_{N+1-i}^{(2N)})^2 - 1}{4} \right )^{j-1} \right ]_{i,j=1,\dots,N}.
$$
The important feature is that each row in the determinant depends on a
distinct coordinate, and allows the sum in (\ref{2.15a}) to be carried
out by summing each $x_{N+1-i}^{(2N)}$ from 1 to $x_{N+1-i}^{(2N+1)} - 1$
over odd values, or equivalently summing $(x_{N+1-i}^{(2N)} - 1)/2$ from
zero to $x_{N+1-i}^{(2N+1)}/2 - 1$ over integer values. Thus
$$
\card
K(x^{(2N)}) = Z_{2N}^{-1} \det \left [ \sum_{x=0}^{x_{N+1-i}^{(2N+1)}/2 - 1}
((x+1)x)^{j-1} \right ]_{i,j=1,\dots,N}.
$$
The sum in column $j$ is an odd polynomial in $x_{N+1-i}^{(2N+1)}$ of degree
$2j-1$ with leading term $(1/(2j-1))(x_{N+1-i}/2)^{2j-1}$ in row $i$. Adding
appropriate multiples of columns $1,\dots,j-1$ to column $j$ shows that this
term can replace the sum, and so
$$
\card K(x^{(2N+1)}) = Z_{2N}^{-1} \prod_{j=1}^N \frac{1}{2j-1} \,
 \det \left [ (x_{N+1-i}/2)^{2j-1} \right ]_{i,j=1,\dots,N}.
$$
Removing a common factor of $x_{N+1-i}/2$ from each row and further use of
the Vandermonde determinant identity gives (\ref{2.16}).

Analogous working shows that with (\ref{2.16}) established in the case $n=2N+1$,
(\ref{2.15}) follows in the case $n=2N+2$.
\end{proof}

\subsection{Asymptotics}
\begin{thm}\label{As}
Under the rescaling $x^{(n)}_i  = \sqrt{2N(1-1/\sqrt{3})} \lambda^{(n)}_i$
the measure $\mathcal{P}^{n}_{\rm blue}$ converges weakly to 
\begin{equation}
\mathcal{P}^{(n)}_{\rm aGUE}\{ \lambda^{(n)}_1,\dots,\lambda^{(n)}_{n/2} \}
= W_n^{-1} 
\prod_{1\leq i<j\leq n/2}
 ((\lambda_j^{(n)})^2 - (\lambda_i^{(n)})^2 )^2
\prod_{ i=1}^{ n/2}
e^{-(\lambda_i^{(n)})^2}
\end{equation}
if $n$ is even and 
\begin{equation}
\mathcal{P}^{(n)}_{\rm aGUE}\{ \lambda^{(n)}_1,\dots,\lambda^{(n)}_{\lfloor n/2\rfloor} \}
= W_n^{-1} 
\prod_{1\leq i<j\leq\lfloor n/2\rfloor}
 ((\lambda_j^{(n)})^2 - (\lambda_i^{(n)})^2 )^2
\prod_{ i=1}^{\lfloor n/2\rfloor}
(\lambda_i^{(n)})^2 e^{-(\lambda_i^{(n)})^2}
\end{equation}
if $n$ is odd. $W_n$ is the normalisation constant $W_n = 
\prod_{l=1}^{\lfloor n/2 \rfloor} N_{n - 2l}$ where $N_j$ is specified by
(\ref{4.37}) below.
\end{thm}
\begin{proof}
Apply Stirling's approximation to the 
formulas from  theorem~\ref{thm:blue}.
\end{proof}

As indicated by the use of the subscript aGUE, $\mathcal{P}^{(n)}_\text{aGUE}$ happens
to be the eigenvalue measure for the anti-symmetric GUE ensemble, 
as we shall see in the next section. 
We are interested in the full measure $\mathcal{P}_\text{blue}^{(1,n)}$,
and how its limit relates to aGUE matrices.

Let $\lambda^{(k)}=(\lambda^{(k)}_1, 
\dots,\lambda^{(k)}_{\lfloor k/2 \rfloor} ) 
\in (\mathbb{R}^+)^n$, for $k=1$, $2$, \dots .
Also let $\lambda^{(2l+1)}_l=0$ for $l=1$, $2$, \dots.
Consider the cone 
 $\mathbf{K}(\lambda^{(n)})=
\{(\lambda^{(1)},\dots, \lambda^{(n)}): \lambda^{(n)}\succ\lambda^{(n-1)}\succ\dots\succ\lambda^{(1)}\}$ 
where $\lambda^{(k+1)}\succ\lambda^{(k)}$ means that 
$\lambda^{(k+1)}_1>\lambda^{(k)}_1>\lambda^{(k+1)}_2>\dots$. 
\begin{thm}
Under the rescaling $x^{(n)}_i = \sqrt{2N(1-1/\sqrt{3})} \lambda^{(n)}_i$, 
as $N=p\rightarrow\infty$, 
the measure $\mathcal{P}^{(n)}_\text{blue}$ converges 
weakly to 
\begin{equation}
\mathcal{P}_{\rm aGUEm}^{(1,n)} \{ \lambda^{(1)},\dots, \lambda^{(n)}\} = 
\frac{\chi_1(\lambda^{(1)}, \lambda^{(2)})
\dots
\chi_{n-1}(\lambda^{{(n-1)}}, \lambda^{(n)})}{
\vol \mathbf{K}(\lambda^{(n)} )}
\mathcal{P}^{(n)}_{\rm aGUE} \{
\lambda^{(n)}
\}
\end{equation}
where $\chi_k(\lambda^{(k)}, \lambda^{(k+1)})=\mathbf{1}\{ \lambda^{(k+1)}\succ \lambda^{(k)} \}$. 
\end{thm}
\begin{proof}
The convergence of a Riemann sum to an integral. 
The number of integer point in a cone $K(x^n)$ suitably normalised
converges to the volume of the cone $\mathbf{K}(\lambda^n)$.
\end{proof}

\begin{lemma}
Let $Z_n$ be as in lemma \ref{lem2.6}. One has
\begin{equation}
\vol \mathbf{K}(\lambda^{(n)}) = Z_n^{-1} \prod_{i<j} ((\lambda^{(n)}_i)^2-(\lambda^{(n)}_j)^2)
\end{equation}
for $n$ even and
\begin{equation}
\vol \mathbf{K}(\lambda^{(n)}) = Z_n^{-1} \prod_{i<j} ((\lambda^{(n)}_i)^2-(\lambda^{(n)}_j)^2)
\prod_{i} (\lambda_i^{(n)})
\end{equation}
for $n$ odd. 
\end{lemma}
\begin{proof}
The proof of lemma \ref{lem2.6}, with summation replaced by integration. 
\end{proof}

The interlacing conditions $\chi_k$ can be written in 
terms of a determinant,
\begin{equation}
\chi_k(\lambda^{(k)}, \lambda^{(k+1)}) = 
\det[  \mathbf{1}\{\lambda^{(k)}_i<\lambda^{(k+1)}_j\} ]_{1\leq i,j\leq M_k}
\end{equation}
(see e.g.~\cite{FR02a}*{lemma 1})
where $M_k=\lfloor (k+1)/2\rfloor$ and 
$\lambda^{(2l+1)}_{l}=0$ for $l=0$, $1$, \dots. 
Also, introduce the standard notation for the Vandermonde determinant, 
\begin{equation}
\Delta(a^2) = 
\prod_{1\leq i<j\leq n} \left( a_i^2 -a_j^2\right)
\end{equation}
for $a\in \mathbb{R}^n$.
This enables us to write the limiting measure on a very
nice form.
\comment{label thm:AS}
\begin{cor}
\label{thm:AS}
\begin{equation}
\mathcal{P}_{\rm aGUEm}^{(1,n)} \{ \lambda^{(1)},\dots, \lambda^{(n)}\} = 
\frac{Z_n}{W_n}\Delta((\lambda^{(n)})^2)
\prod_{k=1}^{n-1}\det[  \mathbf{1}\{\lambda^{(k)}_i<\lambda^{(k+1)}_j\} ]
\prod_{i=1}^{n/2} 
e^{-(\lambda^{(n)}_i)^2} 
\end{equation}
for $n$ even and 
\begin{equation}
\mathcal{P}_{\rm aGUEm}^{(1,n)} \{ \lambda^{(1)},\dots, \lambda^{(n)}\} = 
\frac{Z_n}{W_n} \Delta((\lambda^{(n)})^2)
\prod_{k=1}^{n-1}\det[  \mathbf{1}\{\lambda^{(k)}_i<\lambda^{(k+1)}_j\} ]
\prod_{i=1}^{(n-1)/2} 
(\lambda^{(n)}_i) e^{-(\lambda^{(n)}_i)^2} 
\end{equation}
for $n$ odd.
\end{cor}

In the next section the measures $\mathcal{P}_{\rm aGUEm}^{(1,n)}$ will be
identified with the joint eigenvalue PDF for the minors of
anti-symmetric GUE matrices.

\section{Eigenvalue pdf for the anti-symmetric GUE ensemble}

Consider the anti-symmetric GUE ensemble, which is the probability
measure on purely imaginary Hermitian matrices 
with density $Z^{-1}e^{-\Tr H^2/2}$. Here, $Z$ is a normalisation
constant. Equivalently, form a real
Gaussian matrix with entries chosen independently from 
$N(0, 1/\sqrt{2})$ and set $H=\frac{i}{2} (X-X^T)$.
We seek to find the PDF of the eigenvalues of $H$ 
and its principal minors. For this we adapt workings from
\cites{baryshnikov:gue,FR02b} in which this task is carried out for 
$H$ a GUE matrix. First, 
a technical lemma is required.
\comment{label thm:rational\_polynomial}
\begin{lemma}
\label{thm:rational_polynomial}
Let $0<a_1<\dots<a_n$ be fixed real numbers. 
Let $q_1$, \dots, $q_n$ be i.i.d. $\Exp(1)$ random variables,
specified by  the PDF $e^{-x}$ $(x>0)$. 
Consider the random rational function 
\comment{label eq:1}
\begin{equation}
\label{eq:1}
p(\lambda) = \lambda - \sum_{i=1}^{n} \frac{\lambda q_i }{\lambda^2-a_i^2}. 
\end{equation}
It has  $n$ positive zeros denoted $0<b_1<\dots<b_{n}$, and their PDF is 
\clabel{eq:50}\begin{equation}
\label{eq:50}
2^n\frac{\Delta(b^2) }{\Delta(a^2)}\prod_{i=1}^n b_ie^{-b_i^2  +a_i^2} 
\end{equation} 
supported on $a_1<b_1<a_2<\dots<a_n<b_n$.
\end{lemma}
\begin{proof}
Claim 1: $p$ has at least  $n$ positive roots. 
Since the $q_i$ are all non-negative, 
$p(\lambda)\rightarrow -\infty $ as $\lambda \rightarrow a_i+$
and  $p(\lambda)\rightarrow \infty $ as $\lambda \rightarrow a_i-$.
So there must be a root on each of the intervals $(a_i, a_{i+1})$, 
for $i=1,\dots,n-1$. This also accounts for the inequalities 
$a_1<b_1<a_2<\dots<b_{n-1}<a_n$. Further, since $p(\lambda)\rightarrow \infty $
as $\lambda\rightarrow\infty$, there must be a root $b_n>a_n$. 

Claim 2: $p$ has at most  $n$ positive roots. 
The rational function $p$ has simple poles at $\pm a_i$ for 
$i=1,\dots,n$. Let
\begin{equation}
\label{eq:51}q(\lambda)=\prod_{i=1}^n (\lambda^2-a_i^2), 
\end{equation}
then $pq$ is a polynomial of degree $2n+1$ that has the  same zeros as $p$. 
The zeros already accounted for by the above argument are
are $\pm b_1$, \dots, $\pm b_n$ and $0$, so there can be no more. 

It follows from the above that it is possible to write
\comment{label eq:2}
\begin{equation}\label{eq:2}
p(\lambda)=\frac{\lambda\prod_{i=1}^n (\lambda^2-b_i^2)}{
\prod_{i=1}^n (\lambda^2-a_i^2)}. 
\end{equation}
Comparing the residue at $a_i$ of $p$ in~(\ref{eq:1}) and~(\ref{eq:2}),
an elementary computation gives that
\clabel{eq:32}\begin{equation}
\label{eq:32}
-q_i=\frac{\prod_{j=1}^n (a_i^2-b_j^2)}{
\prod^n_{j=1, j\neq i} (a_i^2-a_j^2)}. 
\end{equation}

The PDF for the variables $\{q_i\}_{i=1}^n$
is  
\comment{label eq:3}
\begin{equation}
\label{eq:3}
\exp(-\sum q_i)
\end{equation}
and we want to change variables to $\{b_i\}_{i=1}^n$.
The Jacobian $J$  for that transformation is,
up to a sign, 
\begin{equation}
J=\det\left[\frac{-2b_j q_i}{a_i^2-b_j^2}  \right] = 
\frac{\prod_{1\leq i<j\leq n} (a_i^2-a_j^2)(b_i^2-b_j^2)}{
\prod_{1\leq i,j\leq n} (a_i^2-b_j^2)}
\prod_i (2b_iq_i),
\end{equation}
where the determinant is evaluated with the Cauchy
double alternant identity. 
Inserting the expression for $q_i$ from~(\ref{eq:32})
simplifies this  to
\comment{label eq:4}
\begin{equation}\label{eq:4}
J=2^n\frac{\prod_{1\leq i<j\leq n} (b_i^2-b_j^2)}{
\prod_{1\leq i<j\leq n} (a_i^2-a_j^2)}\prod_{i=1}^n b_i.
\end{equation} 

By expanding~(\ref{eq:1}) and~(\ref{eq:2}) at infinity
and comparing the $1/\lambda$ coefficient it follows that 
\begin{equation}
\label{eq:52}
-\sum_{i=1}^n q_i = \sum_{i=1}^n a_i^2 -b_i^2  .
\end{equation}
Inserting this in~(\ref{eq:3}) and multiplying by the Jacobian~(\ref{eq:4})
gives the sought form~(\ref{eq:50}).
\end{proof}

\comment{label thm:rational\_polynomial\_2}
\begin{lemma}
\label{thm:rational_polynomial_2}
Let $0<a_1<\dots<a_n$ be fixed real numbers. 
Let $q_1$, \dots, $q_n$ be i.i.d. $\Exp(1)$ 
distributed random variables and let $q_0$ be $\Gamma(1/2,1)$ distributed,
having PDF $(\pi x)^{-1/2} e^{-x}$ for $x>0$.
Consider the random rational function 
\comment{label eq:15}
\begin{equation}
\label{eq:15}
p(\lambda) = \lambda - \frac{q_0}{\lambda}- \sum_{i=1}^{n} \frac{\lambda q_i }{\lambda^2-a_i^2}. 
\end{equation}
It has  $n+1$ positive zeros denoted $0<b_0<\dots<b_{n}$ and their PDF is 
\begin{equation}
\label{eq:56}
\frac{2^{n+1}}{\sqrt{\pi}}
\frac{\Delta(b^2) }{\Delta(a^2)}\prod_{i=0}^{n} e^{-b_i^2 } 
\prod_{i=1}^n \frac{e^{a_i^2} }{a_i}
\end{equation}
supported on 
$0<b_0<a_1<b_1<\dots<a_n<b_{n}$.
\end{lemma}
\begin{proof}
It is convenient to introduce $a_0=0$.
Instead of~(\ref{eq:51}), choose
\begin{equation}
q(\lambda)=\lambda \prod_{i=1}^n (\lambda^2-a_i^2).
\end{equation}
The the proof of lemma~\ref{thm:rational_polynomial} goes through 
virtually unchanged with indices starting from zero instead of from one.
The Jacobian expression from~(\ref{eq:4}) can then be simplified
as
\begin{equation}
J=2^{n+1}
\frac{\prod_{0\leq i<j\leq n} (b_i^2-b_j^2)}{
\prod_{0\leq i<j\leq n} (a_i^2-a_j^2)}\prod_{i=0}^n b_i
=2^{n+1}(-1)^n
\frac{\prod_{0\leq i<j\leq n} (b_i^2-b_j^2) \prod_{i=0}^n b_i }{
\prod_{0\leq i<j\leq n} (a_i^2-a_j^2)\prod_{i=1}^n a_i^2}
\end{equation}
Computing the residue at the origin of  $p$ 
in~(\ref{eq:15}) and~(\ref{eq:2}) gives
\begin{equation}
\label{eq:55}\prod_{i=0}^n b_i^2  = q_0 \prod_{i=1}^n a_i^2 .
\end{equation}
The expression corresponding to~(\ref{eq:52}) is
\begin{equation}
\label{eq:54}
-\sum_{i=0}^n q_i = \sum_{i=1}^n a_i - \sum_{i=0}^n b_i.
\end{equation}
The PDF for the variables $\{q_i\}^n_{i=0}$ is 
\begin{equation}
\label{eq:53}
\frac{\exp\left(-\sum_{i=0}^n q_i \right)}{\sqrt{\pi q_0}}.
\end{equation}
Multiplying this with the Jacobian~(\ref{eq:4}),  inserting (\ref{eq:55})
and (\ref{eq:54}) gives the sought form~(\ref{eq:15}).
\end{proof}
\begin{thm}
\label{thm:aGUE}
Let $H$ be an $n\times n$ matrix from the anti-symmetric GUE ensemble.
Let $H_k$ be the $k\times k$ principal minor of $H$. 
Let $\lambda^{(k)}=(\lambda_1^{(k)},\dots, \lambda_{\lfloor k/2\rfloor }^{(k)})$ 
be the positive eigenvalues of $H_k$, ordered so that 
$\lambda_i^k>\lambda_{i+1}^k$.
Then the joint PDF of $\lambda^{(1)},\dots, \lambda^{(n)}$ is 
precisely that given in corollary~\ref{thm:AS}.
\end{thm}

\begin{proof}
Such a matrix $H$ has the property that if $\lambda$ is an 
eigenvalue of $H$, then so is $-\lambda$. 
Also,  if the size of $H$ is odd, this implies that 
one eigenvalue will be zero. 

The proof is an inductive one. 
A $2\times 2$ matrix from this ensemble is of the form 
$(\begin{smallmatrix} 0&a\\-a&0\end{smallmatrix})$ 
where $a\in N(0,1/\sqrt{2})$. 
Its eigenvalues are $\pm a$, confirming the theorem in the 
case $n=2$. 

First, let $n$ be even. 
Consider an  $n\times n$ matrix $A$ from this ensemble. 
The induction assumption is that its eigenvalue PDF is known.
Consider the  $(n+1)\times(n+1)$ matrix given by bordering $A$, 
\begin{equation}
\left( 
\begin{matrix}
A & w \\ w^* & 0
\end{matrix}
\right)
\end{equation}
Here, $w$ is a column vector of $n$ purely imaginary 
numbers, all $N(0,1/\sqrt{2})$. 
The star means transpose and complex conjugate.

The eigenvectors of $A$ can be paired up in the following way. 
If $v$ is an eigenvector corresponding to eigenvalue $\lambda$
then $\bar v$ is an eigenvector corresponding to eigenvalue $-\lambda$. 
Consider a normalised eigenvector, $|v|=1$.  
Since $v$ and $\bar v$ must be orthogonal to each other,
$|\Real v|^2 = \frac{1}{4}(v+\bar v , v+\bar v) = \frac{1}{2} $,
where $(\cdot,\cdot)$ denotes the inner product.

Let $C=[v_1, \bar v_1, v_2, \dots ]$ be the matrix whose columns are all the 
eigenvectors of $A$. 
Then 
\comment{label eq:14}
\begin{equation}
\label{eq:14}
\left( 
\begin{matrix}
C^* & 0 \\ 0 & 1
\end{matrix}
\right)
\left( 
\begin{matrix}
A & w \\ w^* & 0
\end{matrix}
\right)
\left( 
\begin{matrix}
C & 0 \\ 0 & 1
\end{matrix}
\right)
=
\left( 
\begin{matrix}
D & C^*w \\ w^*C & 0
\end{matrix}
\right)
\end{equation}
where $D$ is a diagonal matrix of the eigenvalues. 
It follows from the above considerations of eigenvectors and 
an elementary calculation that 
$w^* C=(a_1, \bar a_1, a_2,\dots) $ where  each $a_i$
is a complex number, the real and imaginary part of which 
are $N(0,1/\sqrt{2})$. 
Let $p_n(\lambda)$ be the characteristic polynomial of $A$
and say that the eigenvalues of $A$ are $\pm \mu_1$, \dots, $\pm\mu_{n/2}$. 
Of course the eigenvalues of $A$ give the factorisation 
of $p_n$ as
\begin{equation}
p_n(\lambda) = (\mu_1^2-\lambda^2) \dots (\mu_{n/2}^2-\lambda^2).
\end{equation}
Then, it can be shown, say by expanding along the last row
of the RHS of (\ref{eq:14}),
that the characteristic polynomial of that larger matrix is 
such that
\begin{equation}
\frac{p_{n+1}(\lambda)}{p_n(\lambda)}
= \lambda - 
\sum_{i=1}^{n/2} \frac{2a_i \bar a_i \lambda }{\lambda_i^2- \lambda^2}.
\end{equation}
A computation shows that $2a_i\bar a_i$ is $\Exp(1)$ distributed.
So we now need to find the PDF of the zeros of this random rational 
function, which is precisely what is given by  
lemma~\ref{thm:rational_polynomial}. 
Multiplying the expression that the 
induction assumption gives us for $n$ with 
the conditional PDF from lemma~\ref{thm:rational_polynomial}
proves the statement for $H_{n+1}$ when $n$ is even. 

Assume now that  $n$ is odd. 
Do the same construction but the matrix $A$ will now have
one eigenvalue which is zero. 
Performing the same bordering as in~(\ref{eq:14}), 
only this time $w^*C=(a_1, \bar a_1, \dots, a_n,\bar a_n, ib)$. 
where $b$ is $N(0,1/2)$. 
As above the characteristic polynomials for the $n\times n$-matrix
and the $(n+1)\times (n+1)$-matrix are related by
\begin{equation}
\frac{p_{n+1}(\lambda)}{p_n(\lambda)}
= \lambda- \frac{b^2}{\lambda} - \sum_{i=1}^{(n-1)/2} 
\frac{2a_i \bar a_i }{\lambda_i^2- \lambda^2}
\end{equation}
Apply this time lemma~\ref{thm:rational_polynomial_2}
to prove the statement for  $H_{n+1}$ for $n$ odd.
\end{proof}

\section{Correlation functions} 
\subsection{The method of Borodin and collaborators}
Consider $2n+1$ discretisations 
$\cM_1$, \dots, $\cM_{2n+1}$ of the half 
interval $[0,\infty)$ each containing 
the point $0$. 
On each set $\cM_j$ distribute points
 $\{x^j_i\}_{i=1}^{\lfloor2j+1/2\rfloor}$ 
with  $x^{2l-1}_i=0$. 
On the configuration of points $\cX:=\bigcup_{j=1}^{2n+1} \{x^j_i\}$
define a (possibly signed) measure by 
\comment{label eq:5}
\begin{multline}
\label{eq:5}
\frac{1}{C} 
\prod_{l=1}^n
\det\left[ W_{2l-1}(x^{2l-1}_i,x^{2l}_j)\right]_{i,j=1,\dots,l}
\times \det\left[ W_{2l}(x^{2l}_i,x^{2l+1}_j)\right]_{i,j=1,\dots,l}\\
\times \det\left[q_{j-1}(x^{2n+1}_k)\right]_{j,k=1,\dots,n}
\end{multline}
for some functions 
$\{W_i:\cM_i\times \cM_{i+1} \rightarrow\mathbb{R}\}_{i=1,\dots,2n}$,
and $\{q_i:\cM_{2n+1}\rightarrow \rightarrow\mathbb{R}\}_{i=1,\dots, n}$.
After making use of the Vandermonde determinant
expansion one sees that the measures in 
Corollary~\ref{thm:AS} are of this general form. 
One viewpoint is that (\ref{eq:5}) specifies a multicomponent system in
which species $l$ $(l=1,\dots,2n+1)$ consists of $\lfloor (l+1)/2 \rfloor$
particles. But for $l$ odd one of these particles is fixed at the origin,
so the number of mobile particles for species $l$ is
$\lfloor l/2 \rfloor$.
We seek an explicit formula for the correlations of the general 
multi-component system specified by~(\ref{eq:5}).

In fact a generalisation of this very problem has been 
addressed and solved in a recent work of Borodin and 
Ferrari,~\cite{borodin:pushasep}. 
This work extends earlier work of Borodin and collaborators
\cites{borodin:eynard_metha,borodin:fluctuation_tasep}. 
However all these developments are relatively new, so it will
be to the benefit of the reader that we give
an account of the derivation herein rather than just quote the result 
(this is further justified as the proof of 
theorem 4.2 in~\cite{borodin:pushasep} is very brief
and calls for a detailed knowledge of the relevant results 
from~\cites{borodin:eynard_metha,borodin:fluctuation_tasep}

To begin introduce the $|\cM_{2n+1}| \times n$ matrix
$\Psi=[q_{j-1}(x_i)]_{x_i \in \cM_{2n+1},j=1,\dots, n}]$, the 
$|\cM_{2l-1}|\times  |\cM_{2l}|$ matrices
$W_{2l-1}:=[W_{2l-1}(x_j,y_k)]_{x_j\in \cM_{2l-1}, y_k\in\cM_{2l}}$
and the $|\cM_{2l}|\times  |\cM_{2l+1}|$ matrices
$W_{2l}:=[W_{2l}(x_j,y_k)]_{x_j\in \cM_{2l}, y_k\in\cM_{2l+1}}$.
Further introduce the $n\times|\cM_{2l}|$ matrices $E_l$
with entries in row $l$ all ones and entries in all other rows all zeros. 
Finally with $\cM:=\{1,\dots,n\}\cup \cM_1\cup \dots\cup \cM_{2n+1}$
define the $|\cM|\times |\cM|$ matrix $L$ by 
\begin{equation}
L=\begin{bmatrix}
\mathbf{0}&\mathbf{0}&E_1&\mathbf{0}&E_2&\dots &E_n&\mathbf{0}\\
\mathbf{0}&\mathbf{0}&-W_1&\mathbf{0}&\mathbf{0}&\dots&\mathbf{0}&\mathbf{0}\\
\mathbf{0}&\mathbf{0}&\mathbf{0}&-W_2&\mathbf{0}&\dots&\mathbf{0}&\mathbf{0}\\
\mathbf{0}&\mathbf{0}&\mathbf{0}&\mathbf{0}&-W_3&\dots&\mathbf{0}&\mathbf{0}\\
\vdots &&&\ddots\\
\mathbf{0}&\mathbf{0}&\mathbf{0}&\mathbf{0}&\mathbf{0}&\dots&-W_{2n-1}&\mathbf{0}\\
\mathbf{0}&\mathbf{0}&\mathbf{0}&\mathbf{0}&\mathbf{0}&\dots&\mathbf{0}&-W_{2n}\\
\Psi&\mathbf{0}&\mathbf{0}&\mathbf{0}&\mathbf{0}&\dots&\mathbf{0}&\mathbf{0}
\end{bmatrix}
\end{equation}
for block zero matrices $\mathbf{0}$ of appropriate dimension. 

Note that the rows and columns of $L$ are labelled
by the elements of the set $\cM$ in order.
For $Y$ a subset of $\cM$, introduce the notation $L_Y$
to denote the restriction of $L$ to the 
corresponding rows and columns.
Then,  in accordance with the general theory of $L$-ensembles 
(see~\cite{borodin:eynard_metha}) one sees that the measure~(\ref{eq:5})
can be written 
\comment{label eq:6}
\begin{equation}
\label{eq:6}
\frac{\det L_{\{1, \dots,n\}\cup \cX}}{
\det (\mathbf{1}^{*} + L)}
\end{equation}
where $\mathbf{1}^{*}$ is the $\cM\times\cM$ identity matrix with the first 
$n$ ones set to zero. 
The significance of this structure is the general fact that 
the correlation function for particles at 
$Y\in\cM_1\cup\dots \cup\cM_{2n+1}$ is given by~\cite{borodin:eynard_metha},
\begin{equation}
\rho(Y)=\det K_Y, \quad K=\mathbf{1}-(\mathbf{1} + L_{\cM\setminus \{1, \dots, n\}})^{-1}
\end{equation}
for $\mathbf{1}$ the identity of appropriate dimension.
The correlation functions are thus given by a determinant of 
size $|Y|$---the number of particles specified in 
the correlation---and so by definition the measure~(\ref{eq:5})
is a determinantal point process.

Following \cites{borodin:fluctuation_tasep,borodin:pushasep}
we will now proceed to isolate special structures in the functional 
form~(\ref{eq:5}) which when present allow $(\mathbf{1}+L)^{-1}$
in~(\ref{eq:6}) to be computed explicitly. For this purpose,
introduce the matrix 
\begin{equation}
W_{[i,j)}= \begin{cases}
W_i\dots W_{j-1}, & i<j\\
\mathbf{0}, & i\geq j.
\end{cases}
\end{equation}

For $p=1, \dots, 2n$ and $j=0,\dots,n-1$ set
\comment{label eq:7}
\begin{equation}
\label{eq:7}
\Psi^p_{p-2j-2}(x)=(W_{[p, 2n+1)} \Psi )_{x, n-j},\quad x\in \cM_p
\end{equation}
(we use the notation $(A)_{a,b}$ to denote the element in row $a$
column $b$ of the matrix $A$). 
Define functions 
$\{\Phi^p_{p-2l-2}\}_{
\begin{smallmatrix}l=0,1,\dots,\\ p-2l-2\geq 0\end{smallmatrix}}$
constructed from 
\begin{equation}
\linspan\{\{ (E_iW_{[2i,p)})_{2i,x}\}_{i=1, \dots, [(p-1)/2]}
\cup\{\delta_{p, {\rm even}}E_{p/2}\}\}
\end{equation} 
by the orthogonality requirement 
\comment{label eq:8}
\begin{equation}
\label{eq:8}
\sum_{x\in\cM_p} \Phi^p_{p-2j-2}(x)\Psi^p_{p-2k-2}(x) = \delta_{j,k}
\end{equation}
for $j,k=0, 1,\dots$ with $p-2j-2$, $p-2k-2\geq 0$. 
Also set $\Psi_j^{2n+1}(x) = q_j(x)$ $(j=0,\dots,n-1)$.
In terms of this notation, theorem 4.2 in~\cite{borodin:pushasep}
as applied to the measure (\ref{eq:5}) isolates on a special
structure which when satisfied allows the elements of 
$K_Y$ in~(\ref{eq:6}) to be made explicit. 

\comment{label thm:borodin}
\begin{prop}[Theorem 4.2 in~\cite{borodin:pushasep}]
\label{thm:borodin}
Define the $[m/2]\times [m/2]$ matrix $B_m$
by 
\comment{label eq:9}
\begin{equation}
\label{eq:9}
(E_iW_{[2i,m)})_{2i,x} + \delta_{i, m/2} (E_{m/2})_{m, x}
= \sum_{l=1}^{[m/2]}
(B_m)_{i,l} \Phi_{m-2l}^m(x). 
\end{equation}
Suppose $B_m$ is upper triangular. With the set
$Y$ in (\ref{eq:6})
labelled $\{(s_l,y_l)\}$ where 
$y_l\in\cM_{s_l}$ we then have 
\comment{label eq:11}
\begin{equation}
\label{eq:11}
K((r, x), (s, y)) := 
(K_Y)_{(r, x), (s,y)} =
-(W_{[r, s)} )_{x, y} + \sum_{l=1}^{[s/2]} \Psi_{r-2l}^r(x) \Phi_{s-2l}^s(y).
\end{equation}
\end{prop}
\begin{proof}
We follow both the proofs of
\cite{borodin:fluctuation_tasep}*{lemma 3.4}
and \cite{borodin:pushasep}*{theorem 4.2},
and make use of \cite{borodin:eynard_metha}. 
The task is to compute the inverse in~(\ref{eq:6}).

Write $L$ in the structured form 
\begin{equation}
L=\begin{bmatrix}
\mathbf{0} & B\\
C & D-\mathbf{1}
\end{bmatrix}
\end{equation}
where 
\begin{align}
B=&[\mathbf{0}, E_1, \mathbf{0}, E_2,\dots, \mathbf{0},E_n, \mathbf{0}]\\
C=&[\mathbf{0},\dots,\mathbf{0},\Psi]^T.
\end{align}
We know from~\cite{borodin:eynard_metha}*{lemma 1.5}
that
\begin{equation}
K=\mathbf{1} -D^{-1} +D^{-1}CM^{-1} BD^{-1}
\end{equation}
with $M:=BD^{-1} C$ and furthermore 
\begin{equation}
D^{-1} = \mathbf{1}+ [W_{i,j)}]_{i,j=1,\dots,2n+1}.
\end{equation}
{}From this latter formula we compute 
\begin{equation}
D^{-1}C=[W_{[1,2n+1)} \Psi, \dots, W_{[2n,2n+1)}\Psi, \Psi ]^T,
\end{equation}
and we compute that the $m$-th member of the block 
row vector $BD^{-1}$ is equal to 
\begin{equation}
\sum_{k=1}^{[m/2]-1}
E_kW_{[2k,m)} + 
\delta_{m, {\rm even}} E_{m/2}. 
\end{equation}
These formulas in turn tell us that the $(j,m)$
block of $K$ is equal to 
\comment{label eq:10}
\begin{equation}
\label{eq:10}
-W_{[j,m)}+
W_{[j,2n+1)}\Psi M^{-1} (\sum_{k=1}^{[m/2]-1} E_kW_{[2k,m)} + 
\delta_{m, {\rm even}} E_{m/2}),
\end{equation}
where $W_{[j, 2n+1)}$ is to be replaced by $\mathbf{1}$ for $j=2n+1$. 
They tell us too that 
\begin{equation}
M=(\sum_{k=1}^{n-1} E_kW_{[2k, m)})\Psi.
\end{equation}

Define the $(2n+1)\times [m/2]$ matrix $\Phi^m$
specified by 
\begin{equation}
(\Phi^m)_{i,j} =\begin{cases}
\Phi^m_{m-2i} (x^m_j) ,
& 1\leq i\leq [m/2]\\
0, &\text{otherwise.}
\end{cases}
\end{equation}
The orthogonality~(\ref{eq:9}) used in conjunction with~(\ref{eq:10})
shows 
\begin{equation}
\sum_{k=1}^{[m/2]-1}E_k W_{[2k,m)} 
+\delta_{m,{\rm even}}E_{m/2}
= \begin{bmatrix}
B_m &\mathbf{0}\\\mathbf{0}&\mathbf{0}\end{bmatrix}
\Phi^m
\end{equation}
and furthermore $M=B_{2n+1}$.
These formulas, together with the assumption 
that $B_m$ is upper triangular, give that~(\ref{eq:10})
reduces to 
\begin{equation}
-W_{[j,m)} + W_{[j, 2n+1)}\Psi\Phi^m.
\end{equation}

It follows from this and the definition~(\ref{eq:7}) that the element in row 
$(r,x)$ column $(s,y)$ is equal to~(\ref{eq:11}), as required.
\end{proof}

As formulated proposition~\ref{thm:borodin}
applies in the setting that the domains have been discretised. 
Minor modification of the definitions gives that 
(\ref{eq:11}) remains valid in the continuum limit. 
For this purpose, introduce the notation
\begin{equation}
(a*b)(x,y)=\int_0^{\infty} a(x,z) b(z,y)\,dz,
\end{equation}
and note that such an operation is the 
continuum limit of matrix multiplication. 
Then specify
\begin{align}
W_{[i,j)}(x,y)&:=\begin{cases}
(W_i*\dots*W_{j-1})(x,y), &i<j\\
0,&i\geq j
\end{cases}\\
\Psi^p_{p-2j-2}(x)
&=(W_{[p,2n+1)}*q_{n-j-1})(x)
\end{align}
and for $p=1, \dots, 2n$ define functions $\{\Phi^p_{p-2l-2}\}_{
\begin{smallmatrix}l=0,1,\dots,\\ p-2l-2\geq 0\end{smallmatrix}}$
constructed from 
\comment{label eq:12}
\begin{equation}
\label{eq:12}
\linspan\left\{
\left\{ \int_0^\infty 
W_{[2i,p)})(t,x)
\right\}_{i=1, \dots, [(p-1)/2]}
\cup\{\delta_{p, {\rm even}}\}
\right\}
\end{equation} 
by the orthogonality 
requirement 
\comment{label eq:13}
\begin{equation}
\label{eq:13}
\int_0^\infty
\Phi^p_{p-2j-2}(x)\Psi^p_{p-2k-2}(x) \, dx=\delta_{j,k}.
\end{equation}

Consider now the explicit form of these 
quantities for the PDF in corollary~\ref{thm:AS}. 
There, independent of $i$, $W_i(x,y)=\chi_{x<y}$
and so 
\clabel{eq:16}
\begin{equation}
\label{eq:16}
W_{[i,j)}(x,y)=\frac{1}{(j-i-1)!}
\chi_{x<y}(y-x)^{j-i-1}.
\end{equation}
Furthermore, from the Vandermonde determinant identity
\begin{equation}
\prod_{i=1}^n x_i \prod_{1\leq j<k\leq n} (x_k^2-x_j^2)=
\det[x^{2k-1}_j]_{j,k=1,\dots,n} 
\propto \det [p_{2k-1}(x_j)]_{j,k=1,\dots,n}
\end{equation}
valid for any polynomials $p_l(x)$ of degree $l$.
We choose the $p_l(x)$ so that 
\begin{equation}
\int_0^\infty e^{- x^2}p_{2j-1} (x) p_{2k-1}(x) \, dx\propto \delta_{j,k}.
\end{equation}
This is achieved by setting $p_l(x)=H_l(x)$, 
where $H_l(x)$ is the Hermite polynomial. 
Thus $q_i(x)=e^{-x^2} H_{2i+1}(x)$
and so 
\begin{equation}
\Psi^p_{p-2j-2}(x)=\frac{1}{(2n-p)!}
\int_x^\infty e^{-y^2}H_{2n-2j-1} (y) (y-x)^{2n-p} \, dy.
\end{equation}
Making use of the Rodrigues formula 
\clabel{eq:23}
\begin{equation}
\label{eq:23}
H_j(x)=(-1)^je^{x^2}\frac{d^j}{dx^j}e^{-x^2}
\end{equation}
it follows from this that for $j\geq 0$
\clabel{eq:19}
\begin{equation}
\label{eq:19}
\Psi^p_j(x)=e^{-x^2}H_j(x)
\end{equation}
while for $j<0$ 
\clabel{eq:18}
\begin{equation}
\label{eq:18}
\Psi^p_j(x)=\frac{1}{(-j-1)!}\int_x^\infty (y-x)^{-j-1}  e^{-y^2} \, dy.
\end{equation}

It follows from~(\ref{eq:16}) that (\ref{eq:12}) has the explicit form 
\begin{equation}
\linspan
\{
\{x^{p-2i}\}_{i=1,\dots,[(p-1)/2]}\cup \{\delta_{p, {\rm even}}\}\}
\end{equation}
{}From this set we can construct 
\clabel{eq:17}
\begin{equation}
\label{eq:17}
\Phi^p_{p-2i} (x) =  \frac{1}{N_{p-2i} } H_{p-2i} (x)
\end{equation}
for $i=1,2,\dots$,  where 
\begin{equation}\label{4.37}
N_{j}=\int_0^\infty  (H_{j} (x))^2 e^{-x^2} \, dx = \sqrt{\pi} 2^{j-1} j!
\end{equation}
which indeed exhibits the orthogonality (\ref{eq:13}). 
Furthermore, the continuum limit of the change of 
basis formula (\ref{eq:9}) reads 
\begin{equation}
x^{m-2i} + \delta_{i, m/2} = \sum^{[m/2]}_{l=1} (B_m)_{i,l}
\Phi^m_{m-2l} (x)
\end{equation}
and this with $\Phi^m_{m-2l} (x)$ a polynomial of 
degree $m-2l$ as in (\ref{eq:17}) gives that $B_m$ is
an upper triangular matrix. We are therefore justified in 
appling proposition~\ref{thm:borodin}.
\begin{prop}
\label{thm:cont_borodin}
Let $\Psi^p_j(x)$ be specified by~(\ref{eq:19}) and~(\ref{eq:18})
for $j\geq 0$ and $j<0$ respectively. Let $\Phi^p_j(x)=\frac{1}{N_j}H_j(x)$ so
as to be consistent with (\ref{eq:17}). The $r$ point correlation is given by
\clabel{eq:31}
\begin{equation}
\label{eq:31}
\rho_{(r)}((s_j, y_j)_{j=1,\dots,r}) = 
\det[K((s_j,y_j),(s_k,y_k))]_{j,k=1,\dots,r}
\end{equation}
with 
\clabel{eq:22}
\begin{equation}
\label{eq:22}
K((s,x), (t,y)) =
-\frac{1}{(t-s-1)!}\chi_{x<y}(y-x)^{t-s-1} + \sum_{l=1}^{[t/2]} 
\Psi^s_{s-2l}(x)\Phi^t_{t-2l}(y)
\end{equation}
Equivalently, for $s\geq t$ 
\begin{equation}
K((s,x),(t,y))=e^{-x^2} \sum_{l=1}^{[t/2]} \frac{H_{s-2l}(x)H_{t-2l}(y)}{N_{t-2l}}
\end{equation}
while for $s<t$
\clabel{eq:20}
\begin{equation}
\label{eq:20}
K((s,x),(t,y))= - e^{-x^2} 
\sum_{l=-\infty}^{0} \frac{H_{s-2l}(x)H_{t-2l}(y)}{N_{t-2l}}.
\end{equation}
\end{prop}
\begin{proof}
The only remaining task is to derive~(\ref{eq:20}). For $s<t$ consider
\clabel{eq:21}
\begin{equation}
\label{eq:21}
- \frac{1}{(t-s-1)!} \chi_{x<|y|} (\sgn y)^t (|y|-x)^{t-s-1}
\end{equation}
which for $y>0$ is equal to the first term in (\ref{eq:22}). 
For $t$ even (odd) this is an even (odd) function of $y$ and so
can be expanded in the basis $\{H_{2j+\epsilon}(y)\}_{j=0, 1,\dots}$
where $\epsilon=0,1$ for $t$ even, odd, according to 
\begin{equation}
- \frac{1}{(t-s-1)!} \sum_{k=0}^\infty \frac{H_{2k+\epsilon}(y)}{N_{2k+\epsilon}}
\int_x^\infty e^{-u^2}(u-x)^{t-s-1}
H_{2k+\epsilon } (u)\,du,
\end{equation}
valid for $x>0$. Making use of the Rodrigues formula~(\ref{eq:23})
and recalling~(\ref{eq:18}) this can be rewritten
\clabel{eq:24}\begin{equation}
\label{eq:24}
-e^{-x^2 } \sum_{k=-\infty}^{[s/2]}
\frac{H_{s-2k}(x) H_{t-2k} (y)}{N_{t-2k}}
-
\sum_{p=[s/2]+1}^{[t/2]} \Psi^s_{s-2p}(x) \Phi^t_{t-2p} (y).
\end{equation}
Substituting  this for the first term in 
(\ref{eq:22}) and writing in the resulting expression 
\begin{equation}
\sum_{l=1}^{[s/2]} \Psi^s_{s-2p}(x) \Phi^t_{t-2p} (y)
= e^{-x^2 } \sum_{k=1}^{[s/2]}
\frac{H_{s-2k}(x) H_{t-2k} (y)}{N_{t-2k}}
\end{equation}
gives (\ref{eq:20}).
\end{proof}

\subsection{The method of Nagao and Forrester}
In the recent work \cite{forrester:nagao}, the correlations
for a class of multicomponent systems generalizing the GUE minor process, 
in which species $l$ consisted of $l$ particles constrained so 
that they interlace with the particles of species 
$l+1$, were computed in two ways. One was by using the 
generalized formula from \cite{borodin:fluctuation_tasep}*{lemma 3.4},
while the other adapted an approach to multi-component determinantal processes
due to Nagao and Forrester \cite{nagao:correlation_dyson}.
In this section we will show how proposition~\ref{thm:cont_borodin}
can be reclaimed by making use of this latter method. 

With $W_i(x,y)=W(x,y)=\chi_{x<y}$, the first step is to rewrite (\ref{eq:5})
so that it reads 
\clabel{eq:25}
\begin{equation} 
\label{eq:25}
\begin{split}
\frac{1}{C} \prod_{l=1}^n
&\det \begin{bmatrix}
\mathbf{1}_{(n-l)\times (n-l)} & \mathbf{0}_{(n-l)\times l}\\
\mathbf{0}_{l\times (n-l)} &\begin{bmatrix}
[W(x_i^{2l-1} , x_i^{2l}) - \kappa_l(x^{2l-1}_i) ]_{
\begin{smallmatrix}
i=1,\dots,l-1\\
j=1,\dots,l
\end{smallmatrix}}\\
[1]_{k=1,\dots,l}
\end{bmatrix}
\end{bmatrix}\\
\times &\det \begin{bmatrix}
\mathbf{1}_{(n-l)\times (n-l)} & \mathbf{0}_{(n-l)\times l}\\
\mathbf{0}_{l\times (n-l)} &
[W(x_i^{2l} , x_i^{2l+1})  ]_{
i,j=1,\dots,l}
\end{bmatrix}\\
\times & \det[q_{j-1}(x_k^{2n+1})]_{j,k=1,\dots,n}.
\end{split}
\end{equation}
Here $\kappa_l(x)$ is arbitrary as the determinant does not depend on 
$\kappa_l$, and $C$ is a normalization which may vary from 
equation to equation below. 
 Furthermore, $q_{j-1}(x)=e^{-x^2} H_{2j-1}(x)$.
Proceeding as in the derivation of the equality
between (\ref{eq:21}) and~(\ref{eq:24})
we deduce 
\begin{equation}
W(x,y) = \frac{1}{N_0} \int_x^\infty e^{-t^2} \, dt+
e^{-x^2} \sum_{k=1}^\infty \frac{H_{2k-1}(x)H_{2k}(y)}{N_{2k}}
\end{equation}
or alternatively 
\begin{equation}
W(x,y) = 
e^{-x^2} \sum_{k=1}^\infty \frac{H_{2k}(x)H_{2k+1}(y)}{N_{2k+1}}
\end{equation}
both valid for $x,y>0$. We substitute the first of these in~(\ref{eq:25}),
after choosing 
\begin{equation}
\kappa_l(x)=\frac{1}{N_0} \int_x^\infty e^{-t^2} \, dt
\end{equation}
therein, and the second of these in the second of the determinants
in~(\ref{eq:25}).

Introduce now the notation 
\clabel{eq:26}
\begin{equation}
\label{eq:26}
\eta_j(x)=\frac{e^{-x^2/2}}{\gamma_j} H_j(x),\quad \gamma_j:=N^{1/2}_j
\end{equation}
so that $\{\eta_{2j}(x)\}_{j=0, 1,\dots}$ and
$\{\eta_{2j+1}(x)\}_{j=0, 1,\dots}$ each form a set of orthonormal functions 
on $[0,\infty)$. Set
\begin{align}
\phi^{(o)}(x,y)&:=
\sum_{k=1}^\infty 
\frac{\gamma_{2k-1}}{\gamma_{2k}} \eta_{2k-1}(x) \eta_{2k}(y)\\
\phi^{(e)}(x,y)&:=
\sum_{k=0}^\infty 
\frac{\gamma_{2k}}{\gamma_{2k+1}} \eta_{2k}(x) \eta_{2k+1}(y).
\end{align}

{}From the above working we then see that (\ref{eq:25})
can be written 
\clabel{eq:27}\begin{equation} 
\label{eq:27}
\begin{split}
\frac{1}{C} \prod_{l=1}^n
&\det \begin{bmatrix}
\mathbf{1}_{(n-l)\times (n-l)} & \mathbf{0}_{(n-l)\times l}\\
\mathbf{0}_{l\times (n-l)} &\begin{bmatrix}
[1]_{k=1,\dots,l}\\
[\phi^{(o)}(x_i^{2l-1} , x_i^{2l})  ]_{
\begin{smallmatrix}
i=1,\dots,l-1\\
j=1,\dots,l
\end{smallmatrix}}
\end{bmatrix}
\end{bmatrix}\\
\times &\det \begin{bmatrix}
\mathbf{1}_{(n-l)\times (n-l)} & \mathbf{0}_{(n-l)\times l}\\
\mathbf{0}_{l\times (n-l)} &
[\phi^{(e)}(x_i^{2l} , x_i^{2l+1} ]_{i,j=1,\dots,l}
\end{bmatrix}\\
\times & \det[\eta_{j-1}(x_k^{2n+1})]_{j,k=1,\dots,n}
\end{split}
\end{equation}
To  proceed further, set
\begin{equation}
\eta_{j,l}^{(s)} = \begin{cases}
\eta_j(x_l^{(s)}), & j\geq 0, l\geq 1\\
\delta_{j,2l-1}, & \text{otherwise}
\end{cases}
\end{equation}
and use this to define 
\begin{align}
A_{j,l}^{(s,t)} &= \sum_{k=-n}^{-1}
\frac{\gamma_{2k+s}}{\gamma_{2k+t}} 
\eta_{2k+s, j-n+[s/2]}^{(s)}
\eta_{2k+t, l-n+[t/2]}^{(t)} \\
G_{j,l}^{(s,t)} &= \sum_{k=-n}^{\infty}
\frac{\gamma_{2k+s}}{\gamma_{2k+t}} 
\eta_{2k+s, j-n+[s/2]}^{(s)}
\eta_{2k+t, l-n+[t/2]}^{(t)}.
\end{align}
One then sees that 
\begin{equation}
\det[\eta_{2j-1}(x^{2n+1}_l)]_{j,l=1,\dots, l}
\propto
\det[A_{j,l}^{(2n+1,1)}]_{j,l=1,\dots,n}=: \det A^{(2n+1, 1)}
\end{equation}
\begin{equation}
\det\begin{bmatrix} 
[ 1 ]_{k=1, \dots, l}\\
[ \phi^{(o)}(x_i^{2l-1}, x_j^{2l} )]_{
\begin{smallmatrix}
i=1, \dots, l-1\\
j=1, \dots, l
\end{smallmatrix} }
\end{bmatrix}
\propto 
\det[G_{j,k}^{(2l-1,2l)}]_{j,k=1,\dots,n}=: \det G^{(2l-1, 2l)}
\end{equation}
\begin{equation}
\det[ \phi^{(e)}(x_i^{2l}, x_j^{2l+1}) ]_{i,j=1,\dots,l}
\propto \det [G_{j,k}^{(2l,2l+1)}]_{j,k=1,\dots,n} =:
\det G^{(2l,2l+1)}.
\end{equation}
Consequently it is possible to rewrite (\ref{eq:27}) as
\begin{equation}
\frac{1}{C}\det 
\begin{bmatrix}
A^{(2n+1, 1)} &A^{(2n+1, 2)} & A^{(2n+1, 3)} &\dots & A^{(2n+1, 2n+1)} \\
\mathbf{0} & -G^{(1,2)} &-G^{(1,3)} &  \dots & -G^{(1,2n+1)} \\
\mathbf{0} & \mathbf{0} &-G^{(2,3)} &  \dots & -G^{(2,2n+1)} \\
\vdots & \vdots &\vdots & \ddots & \vdots\\
\mathbf{0}&\mathbf{0}&\mathbf{0}&&-G^{(2n, 2n+1)}
\end{bmatrix}
\end{equation}

This has the same structure as \cite{forrester:nagao}*{eq. 4.23}.
Simple manipulation detailed in this latter reference shows that 
this in turn can be rewritten to read
\begin{equation}
\frac{1}{C}
\det [F^{(s,t)}]_{s,t=1, \dots, 2n+1} ,\quad
F^{(s,t)} = \begin{cases}
A^{(s,t)} , & s\geq t\\
B^{(s,t)}, & s<t
\end{cases}
\end{equation}
where $B^{(s,t)}:= A^{(s,t)} - G^{(s,t)} $.
Noting from the definitions of $ A^{(s,t)}$ and 
$G^{(s,t)} $ that for $s\geq t$,
$F^{(s,t)}_{j,l}\propto \delta_{j,l}$ for $j,l\leq n- [s/2]$, while 
for $s<t$, $F^{(s,t)}_{j,l}=0$,
$j\leq n-[s/2]$ or $l\leq n-[t/2]$, this reduces to 
\clabel{eq:29}
\begin{equation}
\label{eq:29}
\frac{1}{C}
\det[f^{(s,t)}]_{s,t=1,\dots,2n+1}
\end{equation}
where $f^{(s,t)}$ is the $s\times t$-matrix with entries 
\clabel{eq:28}
\begin{align}
f_{j,l}^{(s,t)}
&= F^{(s,t)}_{j-[s/2]+n, l-[t/2]+n}\\
\label{eq:28}
&=\begin{cases}
\sum_{k=1}^{[t/2]}
\frac{\gamma^{(s)}_{s-2k}}{\gamma^{(t)}_{t-2k}}
\eta_{s-2k}(x_j^{(s)}) \eta_{t-2k}(x_l^{(t)}),  &s\geq t\\
-\sum_{k=-\infty}^{0}
\frac{\gamma^{(s)}_{s-2k}}{\gamma^{(t)}_{t-2k}}
\eta_{s-2k}(x_j^{(s)}) \eta_{t-2k}(x_l^{(t)}),  &s<t.
\end{cases}
\end{align}
These exhibit the reproducing property 
\begin{equation}
\int_0^{\infty}
f_{j,l}^{(s,t)}
f_{l,m}^{(t,u)}\, dx_l^{(t)} = 
\begin{cases}
f_{j,m}^{(s,u)},& \text{$s\geq t\geq u$ or $s<t<u$}\\
0, &\text{otherwise}
\end{cases}
\end{equation}
which allow the integrations required to reduce~(\ref{eq:29})
down to the $r$-point correlation function
to be performed (see e.g. \cite{forrester:book}), giving
\clabel{eq:30}
\begin{equation}
\label{eq:30}
\rho_{(r)}(\{(s_j,x_j)\}_{j=1,\dots, r})=\det[f_{j,k}^{(s_j,s_k)}]_{j,k=1,\dots, r}
\end{equation}
Substituting (\ref{eq:26}) in~(\ref{eq:28})
shows 
\begin{equation}
f_{j,l}^{(s,t)} = e^{-x_l^2/2} K((s, x_j),( t, x_l))
\end{equation}
and this the result~(\ref{eq:31}) of proposition~\ref{thm:cont_borodin}
is reclaimed.

\section{Scaling limits of the correlations }
In the recent study \cite{forrester:nagao} various scaling limits
of the GUE minor process were computed.
Two involved species which differed by a fixed amount.
In this situation the eigenvalue coordinates were chosen 
to correspond to the neighbourhood of the
 spectrum edge in one case, and the bulk of the spectrum in the other.  
A further scaling limit, 
in which the species differ by $O(N^{2/3})$ and the eigenvalue 
coordinates correspond to the spectrum edge, was also analyzed. 

In the first of these situations, 
the scaled correlation kernel was computed to be the Airy kernel
(\ref{Ksoft})
independent of the particle species. 
In the second situation the correlation kernel for the isotropic
bead process (\ref{eq:33}) was found.
In the final situation, the dynamical extension of the 
Airy kernel (\ref{eq:34})  
was obtained as the scaled correlation kernel. 

The anti-symmetric GUE minor process of the present
work permits the above three scalings. 
We will show below that at the soft edge the limiting correlation kernels
(\ref{eq:32}) and~(\ref{eq:34}) 
are reclaimed in the 
cases of the species differing by $O(1)$ and by $O(n^{2/3})$
respectively. 
In the bulk it will be show that the limiting correlation kernel
(\ref{eq:33}) is reclaimed. 
The antisymmetric GUE minor process also permits the well known 
hard edge scaling \cite{forrester:spectrum_edge}, involving
eigenvalues in the neighbourhood of the origin. 
For eigenvalues from the same species, which is an even label,
the limiting correlation kernel is
\clabel{eq:35}\begin{equation}
\label{eq:35}
K^+(x,y):= 2\int_0^1 \sin(\pi ux)\sin(\pi u y)\, du
\end{equation}
while for odd labelled species it is 
\clabel{eq:36}\begin{equation}
\label{eq:36}
K^-(x,y):= 2\int_0^1 \cos(\pi ux)\cos(\pi u y)\, du.
\end{equation}
In the case of hard edge scaling of the correlation kernel with 
variable species label, 
differing by a finite amount, 
a generalization of the integral representations~(\ref{eq:35})
and (\ref{eq:36}) is found. 
\subsection{Soft edge scaling}
In the $N\times N$ GUE  the soft edge scaling corresponds 
to the change of eigenvalue coordinates 
\clabel{eq:37}\begin{equation}
\label{eq:37}
y_i=\sqrt{2N} + \frac{Y_i}{\sqrt{2}N^{1/6}},
\end{equation}
which has the effect of moving the origin to the 
neighbourhood of the largest eigenvalue and scaling
the distances so the inter-eigenvalue 
spacings in ths neighbourhood are order unity. 
The same soft edge scaling (\ref{eq:37}) applies for 
$(2n+1)\times (2n+1)$ anti-symmetric GUE matrices except that 
$N\mapsto 2n$. 
We consider first this scaling with the species differing from $2n+1$
by a constant,
\clabel{eq:38}\begin{equation}
\label{eq:38}
s_i=2n+1-c_i
\end{equation}
\clabel{thm:soft\_edge}
\begin{prop}
\label{thm:soft_edge}
For the soft edge scaling specified by (\ref{eq:37})
with $N\mapsto 2n$ and
(\ref{eq:38}),
\clabel{eq:39}\begin{equation}
\label{eq:39}
\frac{1}{\sqrt{2} n^{1/6}} 
K((s_j, y_j), (s_l,y_l)) \: \sim \:
\frac{a_n(c_j, Y_j)}{a_n(c_l, Y_l)}K^{\rm soft}(Y_j,Y_l) 
\end{equation}
where $K^{\rm soft}$ is given by (\ref{eq:32}) and 
$a_n(c,Y)=e^{-(2n)^{1/3} Y} (2n)^{-c/2}$. 
Hence we have the pointwise limit 
\begin{equation}
\lim_{n\rightarrow\infty}
\left(\frac{1}{\sqrt{2}n^{1/6}}\right)^r
\rho(\{(s_j, y_j)\}_{j=1, \dots,r} )
=
\det[K^{\rm soft}(Y_j,Y_k)]_{j,k=1,\dots, r}
\end{equation}
independent of the particle species.
\end{prop}
\begin{proof}
Substituting~(\ref{eq:17}) in~(\ref{eq:22})
shows that for $s_j\geq s_l$ ($c_j\leq c_l$)
\clabel{eq:40}\begin{equation}
\label{eq:40}
K((s_j, y_j), (s_l, y_l)) = 
\frac{2e^{-y_j^2}}{\sqrt{\pi}} 
\sum_{k=1}^{[s_j/2]} 
\frac{1}{2^{s_l-2k}(s_l-2k)!}
H_{s_j-2k}(y_j) H_{s_l-2k} (y_l).
\end{equation}
As done for the derivation of the corresponding result in the 
case of the GUE, \cite{forrester:nagao},
we make use of the uniform large $N$
expansion \cite{Ol74} 
\begin{multline}
e^{x^2/2} H_{N-k}(x)=
\pi^{1/4} 2^{(N-k)/2 +1/4} ((N-k)!)^{1/2} N^{-1/12}
\times\\
\times
\left (\Ai\left(X+\frac{k}{N^{1/3}}\right) 
+O(N^{-2/3})\begin{cases}
O(e^{-k/N^{1/3}}), &k\geq 0\\
O(1), & k<0
\end{cases} \right )
\end{multline}
in the summand of (\ref{eq:40})
and so obtain 
\clabel{eq:41}\begin{multline}
\label{eq:41}
K((s_j, y_j), (s_l,y_l))\sim
2 e^{-(2n)^{1/3} (Y_j-Y_k)}
2^{-(c_j-c_l)/2}
\sqrt{2}
(2n)^{-1/6}\times\\
\times
\sum_{k=1}^n
\left( \frac{(2n+1-c_j-2k)!}{(2n+1 -c_l-2k)!}\right)^{1/2}
\Ai(Y_j+\frac{2k}{(2n)^{1/3}})
\Ai(Y_l+\frac{2k}{(2n)^{1/3}}).
\end{multline}
We observe that the leading contribution to 
the sum comes from terms $k=O(n^{1/3})$.
In this regime
\clabel{eq:42}\begin{equation}
\label{eq:42}
\left( \frac{(2n+1-c_j-2k)!}{(2n+1 -c_l-2k)!}\right)^{1/2}
\sim 
(2n)^{(s_l-c_j)/2}
\end{equation}
and so this term can be factored
out of the summand, leaving us with a Riemann
sum approximation to the definite integral (\ref{eq:32})
which implies (\ref{eq:39})

We know from \cite{forrester:nagao} that in the case 
$s_j<s_l$ ($c_j>c_l$) the form (\ref{eq:20})
is not appropriate, due to the resulting Riemann
sum not being convergent. Instead, we make use of (\ref{eq:22}),
breaking the sum over $l$ therein
into the ranges $l\in[1, [s/2]]$, $l\in[[s/2]+1,[t/2]]$ to deduce
that for the given scaling $K((s_j, y_j), (s_l, y_l))$
is to leading order given by the right hand side of 
(\ref{eq:40}), but with the upper terminal replaced by $[s_k/2]$. 
This latter difference does not effect the leading asymptotic form,
and so (\ref{eq:39}) is valid in all cases.
\end{proof}

We turn our attention now to soft edge scaling with the species
separated by $O((2n)^{2/3})$.

\clabel{thm:soft\_edge\_2}
\begin{prop}
\label{thm:soft_edge_2}
Scale $s_i$ according to 
\begin{equation}
s_i=2n-2c_i(2n)^{2/3}
\end{equation}
and scale $y_i$ according to 
\begin{equation}
y_i=(2s_i)^{1/2}+ \frac{Y_i}{\sqrt{2}s_i^{1/6}}.
\end{equation}
For large $n$ and $\alpha_n(c,y):=
e^{-(2n)^{1/3}Y} (4n)^{-c(2n)^{2/3}} 
e^{(2n)^{1/3}c^2} e^{2c^3/3}$,
\clabel{eq:43}\begin{equation}
\label{eq:43}
\frac{1}{\sqrt{2}(2n)^{1/6}}
K((s_j, y_j), (s_l,y_l))\sim \frac{\alpha_n(c_j,Y_j)}{\alpha_n(c_l,Y_l)}
K^{\rm soft} ((c_j, Y_j), (c_l,Y_l))
\end{equation}
where $K^{\rm soft}$ is given by (\ref{eq:32}), and hence 
\begin{equation}
\lim_{n\rightarrow\infty}
\left (\frac{1}{\sqrt{2}(2n)^{1/6}}\right)^r
\rho_{(r)} (\{(s_j, y_j)\}_{j=1,\dots, r})
= \det[K^{\rm soft}((c_j, Y_j), (c_l, y_l))]_{j,l,=1, \dots, r}.
\end{equation}
\end{prop}
\begin{proof}
In the case $c_j\leq c_l$, the formula (\ref{eq:41}) still applies, but now 
with $c_i\mapsto 2c_i (2n)^{2/3}$. This latter point means that (\ref{eq:42}) 
is now inappropriate, following \cite{forrester:nagao}*{eq.~(5.40)} we
use instead the large $s$ expansion
\begin{equation}
\left( \frac{(s-k_j)!}{(s-k_l)!}\right)^{1/2}
\sim
s^{(k_l-k_j)/2} e^{(k_j^2-k_l^2)/4s} e^{(k^3_j-k^3_l)/12s^2}
\end{equation}
with $s=2n$, $k_i=c_i+2k$. A Riemann sum
approximating the first integral in (\ref{eq:34})
is obtained, establishing (\ref{eq:43}) in the case $c_j\leq c_l$.

In the case $c_j>c_l$, it follows from (\ref{eq:20}) that 
(\ref{eq:41}) remains true, but with $k$ now ranging from $-\infty$
to $0$, and the RHS multiplied by $-1$. 
Because the resulting Riemann integral is convergent, 
and is precisely the same as that obtained for $c_j\leq c_l$ 
except that the range is over $(-\infty, 0)$, the second integral 
in (\ref{eq:34}) is obtained, establishing (\ref{eq:43}) 
for $c_j>c_l$. 
\end{proof}

\subsection{Hard edge and bulk}
The hard edge scaling is obtained by change of variables
\clabel{eq:44}\begin{equation}
\label{eq:44}
y_i=\frac{\pi Y_j}{2\sqrt{n}},
\end{equation}
so that the mean particle density in the neighbourhood of the origin is 
of order unity. 
We seek to analyze the correlations with this scaling, 
and the species differing from $2n$ as
specified by (\ref{eq:38}).
Taking $Y_j\rightarrow\infty$ 
in ths expression, with the differences
between the $Y_j$ fixed, the bulk correlations follow as 
a limit of the hard edge correlations. 
\begin{prop}\label{prop5.5}
For the hard edge scaling specified by 
(\ref{eq:44}), (\ref{eq:38}),
\clabel{eq:45}\begin{equation}
\label{eq:45}
\frac{\pi}{2\sqrt{n}} K(s_j, y_j;s_l,y_l) \: \sim \:
2^{c_l-c_j}n^{(c_l-c_j)/2} 
K^{\rm hard} ((s_j, y_j), (s_l,y_l))
\end{equation}
where
\clabel{eq:46}\begin{equation}
\label{eq:46}
K^{\rm hard}
((s,x), (t,y)) = 
\begin{cases}
2\int_0^1 u^{t-s} 
\cos\left(\pi u x -\pi \left(\frac{1-s}{2}\right)\right)
\cos\left(\pi u y -\pi \left(\frac{1-t}{2}\right)\right) \, du, 
&t\geq s\\
-2\int_1^\infty u^{t-s} 
\cos\left(\pi u x -\pi \left(\frac{1-s}{2}\right)\right)
\cos\left(\pi u y -\pi \left(\frac{1-t}{2}\right)\right) \, du, 
&t<s
\end{cases}
\end{equation}
Hence 
\clabel{eq:47}\begin{equation}
\label{eq:47}
\lim_{n\rightarrow\infty}
\left(\frac{\pi}{2\sqrt{n}}\right)^r
\rho(\{ (s_j, y_j)\}_{j=1,\dots,r})
= \det[K^{\rm hard}((c_j, Y_j), (c_l, Y_l))]_{j, l=1, \dots, r}.
\end{equation}
Furthermore, with $s_j$ again given by~(\ref{eq:38}),
the bulk scaling is obtained by the change of variables 
\begin{equation}
y_i=\frac{\pi Y_j}{2\sqrt{n}} + a, \quad 0<a<2\sqrt{n}
\end{equation}
and in this scaling limit
\clabel{eq:48}\begin{equation}
\label{eq:48}
\lim_{n\rightarrow\infty}
\left(\frac{\pi}{2\sqrt{n}}\right)^r
\rho(\{ (s_j, y_j)\}_{j=1,\dots,r})
= \det[K^{\rm bead}((c_j, Y_j), (c_l, Y_l))]_{j, l=1, \dots, r}.
\end{equation}
\end{prop}
\begin{proof}
Consider the hard edge scaling, and suppose $c_l>c_j$.
In (\ref{eq:40}) we substitute for the Hermite polynomials 
the uniform asymptotic expansion \cite{szego:orthogonal}
\begin{equation}
\frac{\Gamma(n/2+1)}{\Gamma(n+1)}
e^{-x^2/2}
H_n(x)=\cos(\sqrt{2n+1} x - n\pi/2)+O(n^{-1/2})
\end{equation}
to deduce 
\clabel{eq:49}
\begin{multline}
\label{eq:49}
K((s_j, y_j), (s_l, y_l))
\sim \frac{2}{\sqrt{\pi}}
\sum_{k=1}^{[s_j/2]}
\frac{1}{2^{s_l-2k} }
\frac{(s_j-2k)!}{ (s_j/2-k)!(s_l/2-k)! }
\times\\
\times 
\cos\left( \pi\sqrt{\frac{n-k}{n}}Y_j -\pi\left(\frac{1-c_j}{2}\right)\right) 
\cos\left( \pi\sqrt{\frac{n-k}{n}}Y_l -\pi\left(\frac{1-c_l}{2}\right)\right) .
\end{multline}
Substituting in this (\ref{eq:38}), where still to be done, 
and scaling the summand by
\begin{equation}
2^{2(n-k)}
(\pi(n-k))^{-1/2}\frac{((n-k)!)^2}{(2n-2k)!}
\end{equation}
(which for large $n-k$ tends to unity) shows 
\begin{multline}
K((s_j, y_j), (s_l, y_l)) \sim
\frac{2}{\pi}
2^{c_l-c_j} n^{(-1+c_l-c_j)/2} 
\sum_{k=1}^n 
\left(\frac{n-k}{n}\right)^{(-1+c_l-c_j)/2}\times
\\
\times\cos\left( \pi\sqrt{\frac{n-k}{n}}Y_j -\pi\left(\frac{1-c_j}{2}\right)\right) 
\cos\left( \pi\sqrt{\frac{n-k}{n}}Y_l -\pi\left(\frac{1-c_l}{2}\right)\right) .
\end{multline}
This is a Riemann sum to the first integral in~(\ref{eq:46}),
implying the leading asymptotics (\ref{eq:45}) in this case $c_l\geq c_j$.

Consider next the hard edge scaling in the case $c_l<c_j$.
As in the analogous stage of the proof of proposition~\ref{thm:soft_edge},
it follows from (\ref{eq:20}) that (\ref{eq:49}) remains true but 
with the summation now over $(-\infty, 0]$, and the RHS multiplied 
by $-1$. This gives a Riemann sum approximation to the second integral 
in~(\ref{eq:46}) and establishes (\ref{eq:45}) in the case $c_l<c_j$.

The result (\ref{eq:48}) can be deduced from 
(\ref{eq:47}) by noting that the sought bulk scaling correlation 
kernel must be the limit $x,y\rightarrow\infty$, $|x-y|$ fixed
of $K^\text{hard}((s,x);(t,y))$. 
A simple rewrite of the resulting integrals gives (\ref{eq:33}).
\end{proof}

\subsubsection{Acknowledgements}
This work was supported by the Australian Research Council and the 
Göran Gustafsson Foundation (KVA).

\end{document}